\documentclass[a4paper]{amsproc}

\usepackage{amssymb}
\usepackage{pstricks}
\usepackage{pstcol}
\usepackage{graphicx}
\usepackage{amsmath}
\usepackage{amsfonts}
\usepackage{latexsym}

\theoremstyle{plain}

\newtheorem{prop}{Proposition}[section]
\newtheorem{lem}{Lemma}[section]
\newtheorem{cor}{Corollary}[section]
\theoremstyle{definition}
\newtheorem{exm}{Example}[section]

\numberwithin{equation}{section}

\title[Elliptic integrals]{The integrals in Gradshteyn and Ryzhik. \\
Part 16: Complete elliptic integrals}

\subjclass[2000]{Primary 33}


\author[S. Boettner]{Stefan Boettner}
\address{Department of Mathematics,
Tulane University, New Orleans, LA 70118}
\email{sboettner@math.tulane.edu}

\author[V. Moll]{Victor H. Moll}
\address{Department of Mathematics,
Tulane University, New Orleans, LA 70118}
\email{vhm@math.tulane.edu}

\address{\hfill{\it Received ??, 
revised ?? }\newline Departamento de Matem\'atica
\newline
Universidad T\'ecnica Federico Santa Mar\'{\i}a
\newline  Casilla 110-V,
\newline Valpara\'{\i}so, Chile}

\thanks{The first author wishes to acknowledge the partial support of  
$\text{NSF-DMS } 0713836$ as a graduate student.  The work of the second 
author was partially supported by the same grant.}

\begin{document}

{\begin{flushleft}\baselineskip9pt\scriptsize {\bf SCIENTIA}\newline
Series A: {\it Mathematical Sciences}, Vol. ?? (2010), ??
\newline Universidad T\'ecnica Federico Santa Mar{\'\i}a
\newline Valpara{\'\i}so, Chile
\newline ISSN 0716-8446
\newline {\copyright\space Universidad T\'ecnica Federico Santa
Mar{\'\i}a\space 2010}
\end{flushleft}}

\vspace{10mm} \setcounter{page}{1} \thispagestyle{empty}

\begin{abstract}
The table of Gradshteyn and Ryzhik contains many entries that are related  
to elliptic integrals. We present a systematic derivation of some of 
them.
\end{abstract}

\maketitle

\newcommand{\nn}{\nonumber}
\newcommand{\ba}{\begin{eqnarray}}
\newcommand{\ea}{\end{eqnarray}}
\newcommand{\ift}{\int_{0}^{\infty}}
\newcommand{\ione}{\int_{0}^{1}}
\newcommand{\ifft}{\int_{- \infty}^{\infty}}
\newcommand{\no}{\noindent}
\newcommand{\realpart}{\mathop{\rm Re}\nolimits}
\newcommand{\imagpart}{\mathop{\rm Im}\nolimits}

\newtheorem{Definition}{\bf Definition}[section]
\newtheorem{Thm}[Definition]{\bf Theorem} 
\newtheorem{Example}[Definition]{\bf Example} 
\newtheorem{Lem}[Definition]{\bf Lemma} 
\newtheorem{Note}[Definition]{\bf Note} 
\newtheorem{Cor}[Definition]{\bf Corollary} 
\newtheorem{Prop}[Definition]{\bf Proposition} 
\newtheorem{Problem}[Definition]{\bf Problem} 
\numberwithin{equation}{section}

\maketitle

\section{Introduction} \label{intro} 
\setcounter{equation}{0}

Elliptic integrals were at the center of analysis at the end of 
$19^{th}$-century. The {\bf complete elliptic integral of the first kind} 
defined by 
\begin{equation}
\mathbf{K}(k) := \int_{0}^{1} \frac{dx}{\sqrt{(1-x^{2})(1-k^{2}x^{2})}}
\end{equation}
\noindent
is a function of the so-called {\bf modulus} $k^{2}$.  The corresponding 
{\bf complete elliptic integral of the second  kind} defined by
\begin{equation}
\mathbf{E}(k) := \int_{0}^{1} \sqrt{\frac{1-k^{2}x^{2}}{1-x^{2}}} dx.
\end{equation}
\noindent
The total collection of complete elliptic integrals contains one more, the 
so-called {\bf complete elliptic integral of the third kind} defined by
\begin{equation}
\mathbf{\Pi}(n,k):= \int_{0}^{1} \frac{dx}{(1-n^{2}x^{2}) 
\sqrt{(1-x^{2})(1-k^{2} x^{2})}}.
\end{equation}

The {\bf complementary integrals} are defined by
\begin{equation}
\mathbf{K'}(k) := \mathbf{K}(k')
\end{equation}
\noindent
where $k' = \sqrt{1-k^{2}}$ is the so-called complementary modulus.

The change of variables $x = \sin t$ yields the trigonometric versions 
\begin{equation}
\mathbf{K}(k) = \int_{0}^{\pi/2} \frac{dt}{\sqrt{1 - k^{2} \sin^{2}t}} 
\text{ and }
\mathbf{E}(k) = \int_{0}^{\pi/2} \sqrt{1 - k^{2} \sin^{2}t} \, dt,
\end{equation}
\noindent
with a similar expression for $\mathbf{\Pi}(n,k)$.

In general, an elliptic integral is one of the form
\begin{equation}
I  := \int_{a}^{b} \frac{P(x) \, dx}{y},
\end{equation}
\noindent
where $y^{2}$ is a cubic or quartic polynomial in $x$.  The integral is 
called {\bf complete} if $a$ and $b$ are roots of $y=0$.  It is clear that 
$\mathbf{K}(k)$ is elliptic. The same is true for $\mathbf{E}(k)$, written in 
the form
\begin{equation}
\mathbf{E}(k) := \int_{0}^{1} \frac{(1-k^{2}x^{2}) \, dx}
{\sqrt{(1-x^{2})(1-k^{2}x^{2})}}. 
\label{expE1}
\end{equation}

\section{Some examples} \label{some-examples} 
\setcounter{equation}{0}

In this section we offer some evaluations from \cite{gr} that follow 
directly from the definitions. Some special values are offered first.
The evaluations of these integrals is facilitated by Legendre's relation
\begin{equation}
\mathbf{K}(k) \mathbf{E'}(k) + \mathbf{K'}(k)\mathbf{E}(k) - 
\mathbf{K}(k)\mathbf{K'}(k) = \frac{\pi}{2}.
\end{equation}
\noindent
The reader will find this identity as Exercise $4$ in section $2.4$ of 
\cite{mckmoll}. 

\begin{exm}
\begin{equation}
\mathbf{K} \left( \sqrt{-1} \right) = \frac{1}{4 \sqrt{2 \pi}} \Gamma^{2} 
\left( \frac{1}{4} \right). 
\end{equation}
\noindent
The proof is direct. The integral is
\begin{equation}
\mathbf{K}(\sqrt{-1}) = \int_{0}^{1} \frac{dx}{\sqrt{1-x^{4}}} = 
\frac{1}{4} \int_{0}^{1} y^{-3/4} (1-y)^{-1/2} \, dy = 
\frac{\Gamma(1/4) \, \Gamma(1/2)}{4 \Gamma(3/4)}. 
\end{equation}
\noindent
The result now follows from the symmetry rule 
\begin{equation}
\Gamma(a) \Gamma(1-a) = \frac{\pi}{\sin \pi a}
\end{equation}
\noindent
for the gamma function and the special value $\Gamma(1/2)  = \sqrt{\pi}$. 
This example appear as entry $\mathbf{3.166.16}$ in \cite{gr}. Entry $3.166.18$ states 
that
\begin{equation}
\int_{0}^{1} \frac{x^{2} \, dx}{\sqrt{1-x^{4}}} = 
\frac{1}{\sqrt{2 \pi}} \Gamma^{2} 
\left( \frac{3}{4} \right). 
\end{equation}
\noindent
The proof consists of a reduction to a special value of the beta function. The
change of variables $t = x^{4}$ gives
\begin{equation}
\int_{0}^{1} \frac{x^{2} \, dx}{\sqrt{1-x^{4}}} = 
\frac{1}{4} \int_{0}^{1} t^{-1/4} (1-t)^{-1/2} \, dt.
\end{equation}
\noindent
This integral is $\tfrac{1}{4} B \left( \tfrac{3}{4}, \tfrac{1}{2} \right)$. 
The simplified result is obtained as above. 

Formula (\ref{expE1}) with $k = \sqrt{-1}$ shows that
\begin{equation}
\mathbf{E}(\sqrt{-1}) = \int_{0}^{1} \frac{1+x^{2}}{\sqrt{1-x^{4}}} \, dx.
\end{equation}
\noindent
The values given above show that
\begin{equation}
\mathbf{E}(\sqrt{-1}) = 
\frac{1}{4 \sqrt{2 \pi}} \left[ \Gamma^{2} \left( \frac{1}{4} \right) + 
4 \Gamma^{2} \left( \frac{3}{4} \right) \right].
\end{equation} 
\end{exm}

\begin{exm}
\begin{equation}
\mathbf{K} 
\left( \frac{1}{\sqrt{2}}  \right) = \frac{1}{4 \sqrt{\pi}} \Gamma^{2} 
\left( \frac{1}{4} \right). 
\end{equation}
\noindent
This appears as entry $\mathbf{8.129.1}$ in \cite{gr}. This value comes from 
the previous example and the identity
\begin{equation}
\mathbf{K}(\sqrt{-1}k) = 
\frac{1}{\sqrt{1+k^{2}}} K \left( \frac{k}{\sqrt{1+k^{2}}}
\right), \label{trans-1}
\end{equation}
\noindent
with $k=1$. The identity (\ref{trans-1}) follows by the change of variables 
$x \mapsto x/\sqrt{1+ k^{2}(1-x^{2})}$ in the left-hand side integral. 
\end{exm}

The values of the modulus $k$ for which 
$\mathbf{K'}/\mathbf{K}$ is the square root of an integer
are of considerable interest. These are called the {\bf singular values}. The 
previous example shows that $1/\sqrt{2}$ is the simplest of them: in this 
case 
\begin{equation}
\frac{\mathbf{K'}}{\mathbf{K}} \left( \frac{1}{\sqrt{2}} \right) = 1.
\end{equation}
\noindent 
A list of the first few values $k_{r}$ for which 
\begin{equation}
\frac{\mathbf{K'}}{\mathbf{K}} \left( k_{r} \right) = \sqrt{r}
\end{equation}
\noindent
is given in \cite{borwein1} and it starts with
\begin{equation}
k_{2} = \sqrt{2}-1, \, k_{3} = \frac{\sqrt{2}(\sqrt{3}-1)}{4}, \, 
k_{4} = 3 - 2 \sqrt{2}, \, 
k_{5} = \frac{1}{2} \left( \sqrt{\sqrt{5}-1} - \sqrt{3 - \sqrt{5}} \right).
\nonumber
\end{equation}

\section{An elementary transformation} \label{sec-elem} 
\setcounter{equation}{0}

Elementary manipulations can be emploted to evaluate certain entries in 
\cite{gr}. For instance, direct integration by parts on the integrals 
defining the functions $\mathbf{K}$ and $\mathbf{E}$ produces 
\begin{equation}
\int_{0}^{1} \frac{x \, \text{ arcsin }x}{\sqrt{(1-k^{2}x^{2})^{3}}} \, dx = 
\frac{1}{k^{2}} \left( \frac{\pi}{2k'} - \mathbf{K}(k) \right)
\end{equation}
\noindent
and 
\begin{equation}
\int_{0}^{1} \frac{x \text{ arcsin } x}{\sqrt{1-k^{2}x^{2}}} \, dx = 
\frac{1}{k^{2}} \left( \mathbf{E}(k) - \frac{\pi}{2} k' \right).
\end{equation}
\noindent
This last evaluation appears as entry ${\mathbf{4.522.4}}$ in \cite{gr}. 

On the other hand, several entries in \cite{gr} may be evaluated also 
by integration by parts choosing the inverse trigonometric term to be 
differentiated. Such procedure gives 

\begin{equation}
\int_{0}^{1} \frac{x \, \text{arccos }x \, dx}{\sqrt{1 - k^{2} x^{2}}} = 
\frac{1}{k^{2}} \left( \frac{\pi}{2} - \mathbf{E}(k) \right),
\end{equation}
\noindent
that appears as entry ${\mathbf{4.522.5}}$, 
\begin{equation}
\int_{0}^{1} \frac{x \, \text{arcsin }x \, dx}{\sqrt{k'^{2} + k^{2} x^{2}}} = 
\frac{1}{k^{2}} \left( \frac{\pi}{2} - \mathbf{E}(k) \right),
\end{equation}
\noindent
that appears as entry ${\mathbf{4.522.6}}$, and finally ${\mathbf{4.522.7}}$:
\begin{equation}
\int_{0}^{1} \frac{x \, \text{arccos }x \, dx}{\sqrt{k'^{2} + k^{2} x^{2}}} = 
\frac{1}{k^{2}} \left( - \frac{\pi}{2}k' + \mathbf{E}(k) \right).
\end{equation}

\smallskip

In this section we derive a different type of 
elementary transformation for integrals and 
use it to obtain the value of some elliptic integrals appearing in \cite{gr}.

\begin{lem}
\label{lemma-odd}
Let $f$ be an odd periodic function of period $a$. Then 
\begin{equation}
\int_{0}^{\infty} \frac{f(x)}{x} \, dx = \frac{\pi}{a} 
\int_{0}^{a/2} \frac{f(x)}{\tan \frac{\pi x}{a}} \, dx.
\end{equation}
\end{lem}
\begin{proof}
The result follows by splitting the integral as 
\begin{eqnarray}
\int_{0}^{\infty} \frac{f(x)}{x} \, dx & = & 
\sum_{k=0}^{\infty} \int_{0}^{a} \frac{f(x)}{x+ka} \, dx \nonumber \\
& = & \sum_{k=0}^{\infty} \int_{0}^{a/2} f(x) 
\left[ \frac{1}{x+ka} - \frac{1}{(k+1)a-x} \right] \, dx\nonumber 
\end{eqnarray}
\noindent
and using the partial fraction decomposition 
\begin{equation}
\tan \frac{\pi b}{2} = \frac{4b}{\pi} \sum_{j=1}^{\infty} 
\frac{1}{(2j-1)^{2} - b^{2}},
\end{equation}
\noindent
given as entry $\mathbf{1.421.1}$ in \cite{gr}. 
\end{proof}

\begin{cor}
\label{coro-1}
Let $f$ be an even function with period $a$. Then 
\begin{equation}
\int_{0}^{\infty} \frac{f(x)}{x} \sin \frac{\pi x}{a} \, dx = 
\frac{\pi}{a} \int_{0}^{a/2} f(x) \, dx. 
\end{equation}
In particular, for $a = \pi$, 
\begin{equation}
\int_{0}^{\infty} \frac{f(x)}{x} \sin x \, dx = 
\int_{0}^{\pi/2} f(x) \, dx. 
\label{aequalpi1}
\end{equation}

\end{cor}
\begin{proof}
Apply the lemma to the function $f(x) \sin \tfrac{\pi x}{a}$ which is odd and 
it has period $2a$. The result follows from the half-angle formula
\begin{equation}
\tan \frac{x}{2} = \frac{\sin x}{1+ \cos x}
\end{equation}
\noindent
and the value
\begin{equation}
\int_{0}^{a} f(x) \cos \frac{\pi x}{a} \, dx = 0.
\end{equation}
\end{proof}

A similar results holds for odd functions. These appear 
as entry ${\mathbf{3.033}}$ in \cite{gr}. 

\begin{cor}
\label{coro-2}
Let $f$ be an odd function with period $a$. Then 
\begin{equation}
\int_{0}^{\infty} \frac{f(x)}{x} \sin \frac{\pi x}{a} \, dx = 
\frac{\pi}{a} \int_{0}^{a/2} f(x) \cos \frac{\pi x}{a} \, dx. 
\end{equation}
In particular, for $a = \pi$, 
\begin{equation}
\int_{0}^{\infty} \frac{f(x)}{x} \sin x \, dx = 
\int_{0}^{\pi/2} f(x) \cos x  \, dx. 
\label{aequalpi}
\end{equation}
\end{cor}

\begin{exm}
\label{sinx}
The function $f(x) \equiv 1$ and $a = \pi$ in Corollary \ref{coro-1} gives the 
classical integral 
\begin{equation}
\int_{0}^{\infty} \frac{\sin x}{x} \, dx = \frac{\pi}{2}.
\end{equation}
\noindent
This is entry $\mathbf{3.721.1}$ 
in \cite{gr}. The reader will find in \cite{hardy9,hardy10} a couple of 
articles by G. H. Hardy with an evaluation of the many proofs of this identity.
These papers are available in volume $5$ of his Complete Works. 
\end{exm}

\begin{exm}
Entry  ${\mathbf{3.842.3}}$ 
of \cite{gr} consists of four evaluations, the first of which 
\begin{equation}
\int_{0}^{\infty} \frac{\sin x}{\sqrt{1-k^{2} \sin^{2}x}} \, \frac{dx}{x} 
= \mathbf{K}(k).
\label{trigo-01}
\end{equation}
\noindent
This follows from Corollary \ref{coro-1} by choosing $a= \pi$ 
and $f(x) = 1/\sqrt{1- k^{2} \sin^{2}x}$. A different proof of this 
evaluation is offered in Section \ref{sec-series} below. 
\end{exm}

\begin{exm}
A second integral appearing in ${\mathbf{3.842.3}}$ is 
\begin{equation}
\int_{0}^{\infty} \frac{\sin x}{\sqrt{1-k^{2} \cos^{2}x}} \, \frac{dx}{x} 
= \mathbf{K}(k).
\label{trigo-03}
\end{equation}
\noindent
also follows from Corollary \ref{coro-1}. This is also true for entry 
${\mathbf{3.841.1}}$ 
\begin{equation}
\int_{0}^{\infty} \sin x \, \sqrt{1- k^{2} \sin^{2}x} \, \frac{dx}{x} = 
\mathbf{E}(k)
\end{equation}
\noindent
and its companion entry ${\mathbf{3.841.2}}$
\begin{equation}
\int_{0}^{\infty} \sin x \, \sqrt{1- k^{2} \cos^{2}x} \, \frac{dx}{x} = 
\mathbf{E}(k).
\end{equation}
\end{exm}

\begin{exm}
The elementary method introduced here may be used to evaluate all integrals 
of the type 
\begin{equation}
I_{m,n}(k) := \int_{0}^{\infty} \frac{\sin^{n}x \, \cos^{m}x}
{\sqrt{1 - k^{2} \sin^{2}x}} \, \frac{dx}{x} 
\end{equation}
\noindent
and the companion family
\begin{equation}
J_{m,n}(k) := \int_{0}^{\infty} \frac{\sin^{n}x \, \cos^{m}x}
{\sqrt{1 - k^{2} \cos^{2}x}} \, \frac{dx}{x}.
\end{equation}
\noindent
All entries in the sections ${\mathbf{3.844}}$ and ${\mathbf{3.846}}$ match 
one of these forms. 
\end{exm}

\begin{exm}
Many other evaluations can be produced by this method. For instance,
\begin{equation}
\int_{0}^{\infty} \frac{\sin x \, \log(1 - k^{2} \sin^{2}x) }
{\sqrt{1- k^{2} \sin^{2}x}} \, \frac{dx}{x}  = 
\int_{0}^{\pi/2} \frac{\log(1 - k^{2} \sin^{2}x) }
{\sqrt{1- k^{2} \sin^{2}x}} \, dx. 
\end{equation}
\noindent
The integral on the left appears as entry ${\mathbf{4.432.1}}$
and the one on the right is entry $\mathbf{4.414.1}$ in \cite{gr}. A proof of 
the identity
\begin{equation}
\int_{0}^{\pi/2} \frac{\log(1 - k^{2} \sin^{2}x) }
{\sqrt{1- k^{2} \sin^{2}x}} \, dx = \mathbf{K}(k) \, \ln k', 
\end{equation}
\noindent
is given in Example \ref{proof414}.
\end{exm}

\section{Some principal value integrals} \label{sec-pv} 
\setcounter{equation}{0}

The method described above can be employed to evaluate some entries of 
\cite{gr} provided the integrals are interpreted as Cauchy principal values.

\begin{exm}
The first example is 
\begin{equation}
\int_{0}^{\infty} \frac{\tan x}{\sqrt{1-k^{2} \sin^{2}x}} \, \frac{dx}{x} = 
\mathbf{K}(k),
\label{form-1}
\end{equation}
\noindent
that appears as one of the four entries in  ${\mathbf{3.842.3}}$ of \cite{gr}. 

Let $I_{1}(k)$ denote the integral and introduce the notation
\begin{equation}
f(x) = \frac{\tan x}{\sqrt{1- k^{2} \sin^{2}x}}.
\end{equation}
Then $f(x)$ is odd and it has period $\pi$. The principal value of the integral 
is given by
\begin{equation}
I_{1}(k) = \lim_{\epsilon \to 0} \sum_{j=0}^{\infty} 
\left( \int_{0}^{\pi/2- \epsilon} \frac{f(x)}{x} \, dx + 
        \int_{\pi/2+ \epsilon}^{\pi} \frac{f(x)}{x+ j \pi} \, dx \right). 
\end{equation}
\noindent
The substitution $y = \pi - x$ in the second integral above produces 
\begin{eqnarray}
I_{1}(k) & = & \lim_{\epsilon \to 0} \sum_{j=0}^{\infty} 
\int_{0}^{\pi/2- \epsilon} \left( \frac{1}{x} + 
\frac{1}{x - (j+1) \pi} \right) f(x) \, dx  \nonumber \\
 & = & \lim_{\epsilon \to 0}
 \int_{0}^{\pi/2- \epsilon} \left( \frac{1}{x} + 
\sum_{j=1}^{\infty} \frac{2x}{x^{2}  - j^{2} \pi^{2}} \right) f(x) \, dx.   
\nonumber 
\end{eqnarray}
The series corresponds to the partial fraction expansion of the cotangent 
function. This completes the evaluation of (\ref{form-1}). The reader will 
note that this proof is very similar to that of Lemma \ref{lemma-odd}.

The value 
\begin{equation}
\int_{0}^{\infty} \frac{\tan x}{\sqrt{1-k^{2} \cos^{2}x}} \, \frac{dx}{x} = 
\mathbf{K}(k),
\label{form-4}
\end{equation}
\noindent
that also appears in ${\mathbf{3.842.3}}$ is established using the same 
type of argument. This completes the evaluation of the integrals in that
entry of \cite{gr}. 
\end{exm}

\begin{exm}
Entry ${\mathbf{3.841.3}}$ of \cite{gr} 
\begin{equation}
\int_{0}^{\infty} \tan x \, \sqrt{1 - k^{2} \sin^{2}x} \, \frac{dx}{x} = 
\mathbf{E}(k)
\end{equation}
\noindent
and its companion ${\mathbf{3.841.4}}$ 
\begin{equation}
\int_{0}^{\infty} \tan x \, \sqrt{1 - k^{2} \cos^{2}x} \, \frac{dx}{x} = 
\mathbf{E}(k)
\end{equation}
\noindent
can be established by the method described in the previous example.
\end{exm}

\section{The hypergeometric connection} \label{hyper} 
\setcounter{equation}{0}

The identites among elliptic integrals often make use of the series 
representations 
\begin{equation}
\mathbf{K}(k) = \frac{\pi}{2} \, {_{2}F_{1}} 
\left[ 
\begin{matrix} \tfrac{1}{2} &  & \tfrac{1}{2} \\ & 1 & \end{matrix} ; \, k^{2} 
\right] = 
\frac{\pi}{2} \sum_{j=0}^{\infty} \frac{\left( \tfrac{1}{2} \right)_{j} 
\, \left( \tfrac{1}{2} \right)_{j}}{j!} \frac{k^{2j}}{j!}, 
\label{series-1}
\end{equation}
\noindent
and
\begin{equation}
\mathbf{E}(k) = \frac{\pi}{2} \, {_{2}F_{1}} 
\left[ 
\begin{matrix} - \tfrac{1}{2} &  & 
\tfrac{1}{2} \\ & 1 & \end{matrix} ; \, k^{2} 
\right] = 
\frac{\pi}{2} \sum_{j=0}^{\infty} \frac{\left( -\tfrac{1}{2} \right)_{j} 
\, \left( \tfrac{1}{2} \right)_{j}}{j!} \frac{k^{2j}}{j!}, 
\label{series-2}
\end{equation}
\noindent
where $_{2}F_{1}$ is the classical hypergeometric function 
\begin{equation}
{_{2}F_{1}} 
\left[ 
\begin{matrix} a  &  & b \\ & c & \end{matrix} ; \,  x
\right] = \sum_{j=0}^{\infty} \frac{(a)_{j} \, (b)_{j}}{(c)_{j} \, j!} x^{j}
\end{equation}
\noindent
and 
\begin{equation}
(a)_{j} = a(a+1)(a+2) \cdots (a+j-1),
\end{equation}
\noindent
is the Pochhammer symbol. The value $(a)_{0}=1$ is adopted. 

\section{Evaluation by series expansions} \label{sec-series} 
\setcounter{equation}{0}

In this section we describe a method to evaluate many of the elliptic 
integrals appearing in \cite{gr}. 

\begin{exm}
The first example is entry ${\mathbf{3.842.3}}$ 
\begin{equation}
\int_{0}^{\infty} \frac{\sin x}{\sqrt{1-k^{2} \sin^{2}x}} \, \frac{dx}{x} 
= \mathbf{K}(k)
\label{series-01}
\end{equation}
\noindent
that has been evaluated in Section \ref{sec-elem}.

Define
\begin{equation}
I_{1}(k^{2}) := 
\int_{0}^{\infty} \frac{\sin x}{\sqrt{1-k^{2} \sin^{2}x}} \, \frac{dx}{x}.
\label{series-02}
\end{equation}
\noindent
To evaluate the integral, let $m = k^{2}$ and expand the integrand in power 
series using 
\begin{equation}
\left( \frac{d}{dm} \right)^{j} \frac{\sin x}{x \, \sqrt{1-m \sin^{2}x}} = 
\left( \frac{1}{2} \right)_{j} \frac{\sin^{2j+1}x}{x} 
(1 - m \sin^{2}x )^{-1/2-j}. 
\end{equation}
\noindent
Therefore,
\begin{equation}
I_{1}(m) = \sum_{j=0}^{\infty} \left( \frac{1}{2} \right)_{j} 
\, \frac{m^{j}}{j!} \, 
\int_{0}^{\infty} \frac{\sin^{2j+1}x}{x} \, dx.
\end{equation}
\noindent
The remaining integral is entry 
${\mathbf{ 3.821.7}}$ in \cite{gr}:
\begin{equation}
\int_{0}^{\infty} \frac{\sin^{2j+1}x}{x} \, dx = \frac{(2j-1)!!}{(2j)!!} 
\frac{\pi}{2}.
\label{38217}
\end{equation}
\noindent
The value of the integral (\ref{series-01}) now follows from 
the series representation of $\mathbf{K}(k)$ given in (\ref{series-1}). 

\medskip

\noindent
{\emph{Proof of}}  (\ref{38217}). Start with
\begin{equation}
\sin^{2j+1}x = 2^{-2j} \sum_{\nu=0}^{j} (-1)^{j-\nu} 
\binom{2j+1}{\nu} \sin(2j-2 \nu +1)x
\end{equation}
\noindent
and the integral in Example \ref{sinx} in the form 
\begin{equation}
\int_{0}^{\infty} \frac{\sin \alpha x}{x} \, dx = \frac{\pi}{2}
\end{equation}
\noindent
for $\alpha > 0$, to obtain
\begin{equation}
\int_{0}^{\infty} \frac{\sin^{2j+1}x}{x} \, dx = 
\frac{\pi}{2^{2j+1}} \sum_{\nu=0}^{j} (-1)^{j-\nu} \binom{2j+1}{\nu}. 
\end{equation}

It follows that
\begin{equation}
I_{1}(m) = \frac{\pi}{2} \, \sum_{j=0}^{\infty} \left( \frac{1}{2} \right)_{j} 
\frac{(-1)^{j}}{2^{2j}} \frac{m^{j}}{j!} \times 
\sum_{\nu=0}^{j} (-1)^{\nu} \binom{2j+1}{\nu}.
\end{equation}

The result now follows from the next lemma. 

\begin{Lem}
Let $j, \, k \in \mathbb{N}$. Then 
\begin{equation}
\sum_{\nu=0}^{k} (-1)^{j} \binom{2j+1}{\nu} = (-1)^{k} \binom{2j}{k}.
\end{equation}
\end{Lem}
\begin{proof}
The proof is by induction on $k$. The case $k=0$ is clear. The induction 
hypothesis is used to produce 
\begin{equation}
\sum_{\nu=0}^{k} (-1)^{\nu} \binom{2j+1}{\nu}  = 
(-1)^{k-1} \binom{2j}{k-1} + (-1)^{k} \binom{2j+1}{k}, 
\end{equation}
\noindent
and an elementary calculation reduces this to $(-1)^{k} \binom{2j}{k}$.  This 
completes the proof of (\ref{38217}).
\end{proof}

\medskip

\noindent
{\emph{Second proof of}} (\ref{38217}): apply the identity (\ref{aequalpi}) to 
the function $f(x) = \sin^{2j}x$ to obtain
\begin{equation}
\int_{0}^{\infty} \frac{\sin^{2j+1}x}{x} \, dx = 
\int_{0}^{\pi/2} \sin^{2j}x \, dx.
\end{equation}
\noindent
This last integral is the classical Wallis' formula given by 
\begin{equation}
\int_{0}^{\pi/2} \sin^{2j}x \, dx = \frac{\pi}{2} \frac{ \left( \tfrac{1}{2} 
\right)_{j}}{j!}. 
\end{equation}
\noindent
The reader will find in \cite{aems-1} information about this formula. 
\end{exm}

\medskip

\begin{exm}
Entry ${\mathbf{ 3.841.1}}$ in \cite{gr} 
\begin{equation}
\int_{0}^{\infty} \sin x \, \sqrt{1-k^{2} \sin^{2}x} \, \frac{dx}{x} = 
\mathbf{E}(k)
\end{equation}
\noindent
is established by the same method employed above. The proof starts with the 
expansion of the integrand using
\begin{equation}
\left( \frac{d}{dm} \right)^{j} \frac{\sin x }{x} \sqrt{1-m \sin^{2} x} 
= \left( - \tfrac{1}{2} \right)_{j} \frac{\sin^{2j+1}x}{x} 
(1 - m \sin^{2}x )^{1/2-j}
\end{equation}
\noindent
and then identify the result with (\ref{series-2}). 
\end{exm}

\begin{exm}
Entry ${\mathbf{3.842.4}}$ in \cite{gr} states that 
\begin{equation}
I_{2}(k):= \int_{0}^{\pi/2} \frac{x \, \sin x \, \cos x}{\sqrt{1-k^{2} 
\sin^{2}x}} \, dx = - \frac{\pi k'}{2k^{2}} + \frac{E(k)}{k^{2}}. 
\end{equation}
\noindent
The parameter $k'$ is the complementary modulus $k' = \sqrt{1-k^{2}}$. 

\smallskip

Write $m=k^{2}$ and expand the integrand in series using
\begin{equation}
\left(\frac{d}{dm} \right)^{j} 
\frac{x \sin x \, \cos x}{\sqrt{1-m \sin^{2}x}} = 
\left( \frac{1}{2} \right)_{j} 
\frac{x \sin^{2j+1} x \, \cos x}{\sqrt{1-m \sin^{2}x}}. 
\end{equation}
\noindent
Therefore
\begin{equation}
I_{2}(m) = \sum_{j=0}^{\infty} \left( \frac{1}{2} \right)_{j} 
\frac{m^{j}}{j!} 
\int_{0}^{\pi/2} x \sin^{2j+1}x \, \cos x \, dx.
\end{equation}

Integration by parts gives
\begin{equation}
\int_{0}^{\pi/2} x \sin^{2j+1}x \, \cos x \, dx = 
\frac{\pi}{4(j+1)} - \frac{1}{4(j+1)} B \left( j + \tfrac{3}{2}, \tfrac{1}{2} 
\right),
\end{equation}
\noindent
where 
\begin{equation}
B(u,v) = \int_{0}^{1} t^{u-1}(1-t)^{v-1} \, dt = 
2 \int_{0}^{\pi/2} \sin^{2u-1} \varphi \, \cos^{2v-1} \varphi \, d \varphi,
\end{equation}
\noindent
is the classical beta function.  It follows that
\begin{equation}
I_{2}(m) = \frac{\pi}{4} \sum_{j=0}^{\infty} \left( \frac{1}{2} \right)_{j}
\frac{m^{j}}{(j+1)!} - \frac{1}{4} \sum_{j=0}^{\infty} 
\left( \frac{1}{2} \right)_{j} 
B \left( j + \frac{3}{2}, \frac{1}{2} \right)
\frac{m^{j}}{(j+1)!}. 
\label{formula-i2}
\end{equation}
\noindent
The two series are now treated separately. 

The first sum is computed by the binomial theorem
\begin{equation}
(1-x)^{-a} = \sum_{j=0}^{\infty} \frac{(a)_{j}}{j!} x^{j}
\end{equation}
\noindent
as 
\begin{equation}
\frac{\pi}{4} \sum_{j=0}^{\infty} \left( \frac{1}{2} \right)_{j}
\frac{m^{j}}{(j+1)!} = 
\frac{\pi}{2(1+ \sqrt{1-m})} = \frac{\pi}{2m}(1- \sqrt{1-m}). 
\end{equation}

\smallskip

The second sum is 
\begin{eqnarray}
- \frac{1}{4} \sum_{j=0}^{\infty} 
\left( \frac{1}{2} \right)_{j} 
B \left( j + \frac{3}{2}, \frac{1}{2} \right)
\frac{m^{j}}{(j+1)!} & = & 
- \frac{\sqrt{\pi}}{4} \sum_{j=0}^{\infty} 
\left( \frac{1}{2} \right)_{j} \frac{\Gamma(j + \tfrac{3}{2})}
{(j+1)! \Gamma(j+2)} m^{j} \nonumber \\
& = & - \frac{\pi}{8} \sum_{j=0}^{\infty} 
\left( \frac{1}{2} \right)_{j} 
\left( \frac{1}{2} \right)_{j} 
\frac{m^{j}}{(j+1)!} \nonumber \\
& = &  \frac{\pi}{2m} \sum_{j=0}^{\infty} 
\left( - \frac{1}{2} \right)_{j+1} 
\left( \frac{1}{2} \right)_{j+1} 
\frac{m^{j+1}}{(j+1)!} \nonumber \\
& = & \frac{\pi}{2m} \left[ {_{2}F_{1}} 
\left( \begin{matrix}  -\tfrac{1}{2} & & 
\tfrac{1}{2} \\ & 1 & \end{matrix}; m  \right) -1 \right]. \nonumber
\end{eqnarray}
\noindent
The hypergeometric representation (\ref{series-2}) and (\ref{formula-i2}) give
\begin{equation}
I_{2}(m) = -\frac{\pi \sqrt{1-m}}{2m} + \frac{\mathbf{E}(k)}{m}
\end{equation}
\noindent
as claimed. 
\end{exm}

\section{A small correction to a formula in Gradshteyn and Ryzhik} 
\label{sec-error} 
\setcounter{equation}{0}

In this section we present the evaluation of some elliptic integrals in 
\cite{gr}. In particular, a small error in formula ${\mathbf{4.395.1}}$ is 
corrected.

\begin{prop}
\label{correct}
Let $k' = \sqrt{1-k^{2}}$ be the complementary modulus. Then 
\begin{equation}
\int_{0}^{\infty} \frac{\ln x \, dx}{\sqrt{(1+x^{2})(k'^{2} + x^{2})}}  
= \frac{1}{2} \mathbf{K}(k) \, \ln k'.
\label{formula-22}
\end{equation}
\end{prop}
\begin{proof}
Let $m = k'^{2}$ and use
\begin{equation}
\left( \frac{d}{dm} \right)^{j} \frac{\ln x}{\sqrt{(1+x^{2})(m+x^{2})}} = 
(-1)^{j} \left( \frac{1}{2} \right)_{j} \, 
\frac{\ln x}{\sqrt{(1+x^{2})(m+x^{2})^{j+1/2}}} 
\end{equation}
\noindent
to expand the integrand around $m=-1$. It follows that
\begin{equation}
\int_{0}^{\infty} \frac{\ln x \, dx}{\sqrt{(1+x^{2})(k'^{2} + x^{2})}}   = 
\sum_{j=0}^{\infty} \frac{(-1)^{j}}{j!} \left( \tfrac{1}{2} \right)_{j} 
\int_{0}^{\infty} \frac{\ln x \, dx}{(1+x^{2})^{j+1}} \, 
(m-1)^{j}.
\end{equation}
\noindent
This last integral is given by
\begin{eqnarray}
\int_{0}^{\infty} \frac{\ln x \, dx}{(1+x^{2})^{j+1}} & = & 
\frac{1}{4} \int_{0}^{\infty} \frac{\ln x \, dx}{\sqrt{x} \, (1+ x)^{j+1}} 
\nonumber \\
& = & \frac{1}{4} \frac{d}{d \alpha} B(\alpha, j - \alpha+1) \Big{|}_{\alpha = 1/2} \nonumber \\
& = & \frac{1}{4} B \left( \tfrac{1}{2}, j + \tfrac{1}{2} \right) 
\left[ \psi \left( \tfrac{1}{2} \right) - \psi \left( j + \tfrac{1}{2} \right)
\right] \nonumber \\
& = & \frac{\pi}{2 j!} \left( \frac{1}{2} \right)_{j} 
\, \sum_{i=0}^{j-1} \frac{1}{2i+1}. \nonumber
\end{eqnarray}
\noindent
Therefore, the left-hand side of (\ref{formula-22}) satisfies
\begin{equation}
\int_{0}^{\infty} \frac{\ln x \, dx}{\sqrt{(1+x^{2})(k'^{2} + x^{2})}}  
= \frac{\pi}{2} \sum_{j=0}^{\infty} 
\frac{\left( \tfrac{1}{2} \right)_{j}^{2}}{j!^{2}}  \sum_{i=0}^{j-1} \frac{1}{2i+1} 
(1 - m )^{j}.
\end{equation}
\noindent 
The series expansion for the complete elliptic integral now shows that the 
right-hand side of (\ref{formula-22}) is given by
\begin{eqnarray}
\frac{1}{4} \ln m \, {\mathbf{K}}(\sqrt{1-m}) & = & 
\frac{\pi}{8} \left[ \sum_{j=1}^{\infty} \frac{(1-m)^{j}}{j} \right]
\times  \left[ \sum_{j=0}^{\infty} 
\frac{\left( \tfrac{1}{2} \right)_{j}^{2} }{j!^{2}}
(1-m)^{j} \right] \nonumber \\
& = & \frac{\pi}{8} \sum_{j=0}^{\infty} 
\left[ \sum_{i=0}^{j-1}  \frac{1}{j-i} 
\frac{\left( \tfrac{1}{2} 
\right)_{i}^{2}}{i!^{2}}  \right] (1-m)^{j}.
\nonumber
\end{eqnarray}

The result follows from the identity established in the next lemma.
\end{proof}

\begin{lem}
\label{harmonic}
Let $j \in \mathbb{N}$. Define 
\begin{equation}
a_{r} = \frac{\left( \tfrac{1}{2} \right)_{r}^{2}}{r!^{2}}.
\end{equation}
\noindent
Then 
\begin{equation}
\sum_{i=0}^{j-1} \frac{a_{i}}{j-i} 
= 4a_{j} \,  \sum_{i=0}^{j-1} \frac{1}{2i+1}. 
\end{equation}
\end{lem}
\begin{proof}
The relations 
\begin{equation}
(-x)_{k} = (-1)^{k} (x-k+1)_{k} \text{ and }
\left( \tfrac{1}{2} \right)_{n-k} 
\left( \tfrac{1}{2} - n \right)_{k}  = (-1)^{k} 
\left( \tfrac{1}{2} \right)_{n} 
\end{equation}
\noindent
can be used to rewrite the left-hand side as 
\begin{eqnarray}
\sum_{i=0}^{j-1} \frac{ \left( \tfrac{1}{2} \right)_{i}^{2} }{i!^{2}} \, 
\frac{1}{j-i} & = &  
\sum_{k=0}^{j-1} \frac{ \left( \tfrac{1}{2} \right)_{j-k-1}^{2} }
{(j-k-1!^{2}} \, 
\frac{1}{k+1}  \nonumber \\
& = &  \frac{ \left( \tfrac{1}{2} \right)_{j}^{2} }
{j!^{2}} \, \sum_{k=0}^{j-1} 
\frac{ (-j)_{k+1}^{2}}{ \left( \tfrac{1}{2} - j \right)_{k+1}^{2}} \, 
\frac{1}{k+1}. \nonumber
\end{eqnarray}
\noindent
Thus the assertion of the lemma is equivalent to 
\begin{equation}
\sum_{k=0}^{j-1} \frac{(-j)_{k+1}^{2}}{ \left( \tfrac{1}{2} - j 
\right)_{k+1}^{2}} \, \frac{1}{k+1} = \sum_{k=0}^{j-1} \frac{4}{2k+1}.
\end{equation}

\smallskip

Next apply the fact that $(x)_{k+1} = x (x+1)_{k}$ to obtain

\begin{eqnarray}
\sum_{k=0}^{j-1} \frac{(-j)_{k+1}^{2}}{ \left( \tfrac{1}{2} - 
j \right)_{k+1}^{2}} \, \frac{1}{k+1} & = & 
\frac{j^{2}}{ \left( \tfrac{1}{2} - j \right)^{2}} 
\sum_{k=0}^{j-1} \frac{(1-j)_{k}^{2}}{ \left( \tfrac{3}{2} - 
j \right)_{k}^{2}} \, \frac{1}{k+1}  \nonumber \\
& = & \frac{j^{2}}{ \left( \tfrac{1}{2} - j \right)^{2}} 
\sum_{k=0}^{j-1} \frac{(1-j)_{k}^{2} (1)_{k}^{2}}{ \left( \tfrac{3}{2} - 
j \right)_{k}^{2} (2)_{k} \, k!}.   \nonumber 
\end{eqnarray}

The right-hand side is a balanced $_{4}F_{3}$ series and it can be transformed 
using 
\begin{multline}
{_{4}F_{3}} 
\left[ 
\begin{matrix} x &  & y & & z & & -m  \\ & u & & v & & w \end{matrix}; 1
\right] =  \\
\frac{(v-z)_{m} \, (w-z)_{m}}{(v)_{m}\, (w)_{m}} \, 
{_{4}F_{3}} 
\left[ 
\begin{matrix} u-x &  & u-y & & z & & -m  \\ & 1-v+z-m & & 1-w+z-m 
& & u \end{matrix}; 1
\right]. 
\nonumber
\end{multline}
%The argument in the power series representation of the ${_{4}F_{3}}$ is not 
%written and it is taken at value $1$. 
See \cite{bailey35}, page 56. Now let $y=z=1, \, x = 1-j, \, m = j-1, \, 
u = v = \tfrac{3}{2} -j$ and $w=2$. It follows that
\begin{multline}
{_{4}F_{3}} 
\left[ 
\begin{matrix} 1 &  & 1 & & 1-j & & 1-j  \\ & \tfrac{3}{2}-j & 
& \tfrac{3}{2} -j  & & 2 \end{matrix} \; ; 1
\right]   =  \nonumber \\  
\frac{(\tfrac{1}{2} -j)_{j-1} \, (1)_{j-1}}{(\tfrac{3}{2} - j)_{j-1}\, 
(2)_{j-1}} \, 
{_{4}F_{3}} 
\left[ 
\begin{matrix} 1 &  & \tfrac{1}{2} & & -j + \tfrac{1}{2} & & 1-j  
\\ & -j + \tfrac{3}{2} & &  \tfrac{3}{2}
& & 1-j \end{matrix} \, ; 1
\right]. 
\nonumber
\end{multline}

\noindent
The last hypergeometric terms is now simplified
\begin{eqnarray}
\frac{2j-1}{j} \sum_{k=0}^{j-1} \frac{ \left( \tfrac{1}{2} \right)_{k} \,
\left( -j + \tfrac{1}{2} \right)_{k} }
{ \left( \tfrac{3}{2} \right)_{k} \, \left( - j + \tfrac{3}{2} 
\right)_{k}} & =  & 
\frac{(2j-1)^{2}}{j} \sum_{k=0}^{j-1} \frac{1}{(2k+1)(2j-1-2k)} \nonumber \\
& = & \frac{(2j-1)^{2}}{2j^{2}} \sum_{k=0}^{j-1} 
\left( \frac{1}{2k+1} + \frac{1}{2j-1-2k} \right) \nonumber \\
& = & \frac{(2j-1)^{2}}{j^{2}} \sum_{k=0}^{j-1} \frac{1}{2k+1}, \nonumber 
\end{eqnarray}
\noindent
as claimed. 
\end{proof}

\noindent
{\emph{An automatic proof}.} The result of Lemma \ref{harmonic} also admits an 
automatic proof as described in \cite{aequalsb}. Define the functions $F(i,j)$
and $G(i,j)$, respectively, as 
\begin{equation}
F(i,j) = \frac{\left( \tfrac{1}{2} \right)_{i}^{2} \, j!^{2}}
{\left( \tfrac{1}{2} \right)_{j}^{2} \, i!^{2}} \, \frac{1}{j-i} 
\text{ and }
G(i,j) = - \frac{\left( \tfrac{1}{2} \right)_{i}^{2} \, j!^{2}}
{\left( \tfrac{1}{2} \right)_{j+1}^{2} \, i!^{2}} \, \frac{i^{2}}{j-i+1}.
\end{equation}
\noindent
The stated result is equivalent to the identity $a(j) = b(j)$, where
\begin{equation}
a(j) = \sum_{i=0}^{j-1} F(i,j) \text{ and } b(j) = \sum_{i=0}^{j-1} 
\frac{1}{2j+1}.
\end{equation}
\noindent
Zeilberger algorithm finds the non-homogeneous recurrence 
\begin{equation}
F(i+1,j) - F(i,j) = G(i+1,j)-G(i,j).
\end{equation}
\noindent
Summing this for $i$ from $0$ to $j-1$ and using the telescoping of the 
right-hand side, produces
\begin{eqnarray}
\sum_{i=0}^{j-1} F(i,j+1) - \sum_{i=0}^{j-1} F(i,j) & = & 
\sum_{i=0}^{j-1} G(i+1,j) - \sum_{i=0}^{j-1} G(i,j)  \nonumber \\
& = & G(j,j) - G(0,j) \nonumber \\
& = & - \frac{4j^{2}}{(2j+1)^{2}}. \nonumber
\end{eqnarray}
\noindent
Now observe that
\begin{eqnarray}
a(j+1)-a(j) & = & \frac{4(j+1)^{2}}{(2j+1)^{2}} + 
\sum_{i=0}^{j-1} F(i,j+1) - \sum_{i=0}^{j-1} F(i,j)  \nonumber \\
& = & \frac{4(j+1)^{2}}{(2j+1)^{2}} - \frac{4j^{2}}{(2j+1)^{2}} 
\nonumber \\
 & = & \frac{4}{2j+1}. \nonumber
\end{eqnarray}
\noindent
The sequence $b(j)$ satisfies the same recurrence.
Therefore $a(j)-b(j)$ is a constant. Since $a(1)=b(1) = 4$ 
this constant vanishes. This establishes the result.

\medskip

The next result corrects entry ${\mathbf{4.395.1}}$ in \cite{gr}. 

\begin{cor}
The value 
\begin{equation}
\int_{0}^{\infty} \frac{\ln \tan \theta \, d \theta}
{\sqrt{1- k^{2} \sin^{2} \theta}} = 
- \frac{1}{2} \ln k' \mathbf{K}(k)
\end{equation}
\noindent
holds. 
\end{cor}
\begin{proof}
Let $x \mapsto \tan \theta$ in (\ref{formula-22}).
\end{proof}

\begin{exm}
Entry ${\mathbf{4.242.1}}$ states 
\begin{equation}
\int_{0}^{\infty} \frac{\ln x \, dx }{\sqrt{(a^{2}+x^{2})(b^{2}+x^{2})}} =
\frac{1}{2a} \mathbf{K} \left( \frac{\sqrt{a^{2}-b^{2}}}{a} \right) \, \ln ab.
\end{equation}
\noindent
Formula (\ref{formula-22}) corresponds to the special case $a=1$. The change 
of variables $x = at$ produces
\begin{equation}
\int_{0}^{\infty} \frac{\ln x \, dx }{\sqrt{(a^{2}+x^{2})(b^{2}+x^{2})}} =
\frac{1}{b} \int_{0}^{\infty} \frac{\ln t \, dt}
{\sqrt{(1+t^{2})(c^{2}+t^{2})}} + \frac{\ln a}{b} 
\int_{0}^{\infty} \frac{dt}{\sqrt{(1+t^{2})(1+c^{2}t^{2})}}
\nonumber
\end{equation}
\noindent
with $c = b/a$. The first integral is evaluated using (\ref{formula-22}) 
and let $t = \tan \varphi$ to see that the second integral is
${\mathbf{K}}(\sqrt{1-c^{2}})$. This establishes the result.
\end{exm}

\begin{exm}
\label{proof414}
The techniques illustrated here are now employed to prove entry 
${\mathbf{4.414.1}}$ in \cite{gr}:
\begin{equation}
\int_{0}^{\pi/2} \frac{\ln(1 - k^{2} \sin^{2}x)}{\sqrt{1- k^{2} \sin^{2}x}}
\, dx = {\mathbf{K}}(k) \ln k'.
\label{formula-23}
\end{equation}
\noindent
Let $m = k^{2}$ and observe that 
\begin{equation}
\frac{d}{dm} \frac{\alpha_{j} + \beta_{j} \ln(1 - m \sin^{2}x)}
{(1 - m \sin^{2}x)^{j+1/2}} \sin^{2j}x = 
\frac{\alpha_{j+1} + \beta_{j+1} \ln(1 - m \sin^{2}x)}
{(1 - m \sin^{2}x)^{j+3/2}} \sin^{2j+2}x 
\end{equation}
\noindent
where the parameters $\alpha_{j}, \, \beta_{j}$ satisfy 
\begin{equation}
\alpha_{j+1} = (j + \tfrac{1}{2}) \alpha-{j} - \beta_{j} \text{ and }
\beta_{j+1} = ( j + \tfrac{1}{2} ) \beta_{j}.
\end{equation}
\noindent
Now choose $\alpha_{0}=0$ and $\beta_{0}=1$ to obtain
\begin{equation}
\left( \frac{d}{dm} \right)^{j} \frac{\ln(1 - m \sin^{2}x)}
{\sqrt{1- m \sin^{2}x}} = 
\frac{\alpha_{j} + \beta_{j} \ln(1 - m \sin^{2}x)}{(1-m \sin^{2}x)^{j+1/2}} 
\sin^{2j}x.
\end{equation}
\noindent
Expand the integrand of (\ref{formula-23}) around $m=0$ and use
\begin{equation}
\int_{0}^{\infty} \sin^{2j}x \, dx = \frac{\pi}{2} \frac{\left( \tfrac{1}{2} 
\right)_{j}}{j!^{2}}
\end{equation}
\noindent
and the expressions 
\begin{equation}
\alpha_{j} = \left( \tfrac{1}{2} \right)_{j} \sum_{i=0}^{j-1} \frac{2}{2i+1}
\text{ and } \beta_{j} = \left( \tfrac{1}{2} \right)_{j}
\end{equation}
\noindent
to see that 
\begin{equation}
\int_{0}^{\pi/2} \frac{\ln(1 - k^{2} \sin^{2}x)}{\sqrt{1- k^{2} \sin^{2}x}}
\, dx = 
\pi \sum_{j=0}^{\infty} \frac{ \left( \tfrac{1}{2} \right)_{j}^{2}}{j!^{2}} 
\, \left( \sum_{i=0}^{j-1} \frac{1}{2i+1} \right) \, m^{j}. 
\end{equation}
\noindent
The result now follows from the evaluation given in the proof of Proposition
\ref{correct}. 
\end{exm}

\bigskip

\noindent
{\bf Acknowledgments}. The authors wish to thank T. Amdeberhan for many 
comments on a first draft of this paper. In particular, both proofs of 
Lemma \ref{harmonic} are due to him. The first author acknowledges the 
partial support of 
$\text{NSF-DMS } 0713836$ as a graduate student.  The work of the second 
author was partially supported by the same grant.

\end{document}